\documentclass[11pt]{amsart}
\usepackage{amsfonts}
\usepackage{amsmath}
\usepackage{amssymb}
\usepackage{color}
\usepackage{graphicx}
\usepackage{xspace}
 \font \eightrm=cmr8

 \newcommand{\nc}{\newcommand}

\newtheorem{thm}{Theorem}
\newtheorem{exam}{Example}

\newtheorem{cor}[thm]{Corollary}

\newtheorem{lem}[thm]{Lemma}
\newtheorem{prop}[thm]{Proposition}
\newtheorem{defn}{Definition}
\newtheorem{rmk}[thm]{Remark}

\def\BCH{{\rm{BCH}}}

\nc{\ignore}[1]{{}}
\nc{\mrm}[1]{{\rm #1}}
\nc{\dirlim}{\displaystyle{\lim_{\longrightarrow}}\,}
\nc{\invlim}{\displaystyle{\lim_{\longleftarrow}}\,}
\nc{\vep}{\varepsilon} \nc{\ep}{\epsilon}
\nc{\sigmat}{\widetilde\sigma}
\nc{\ostar}{\overline{*}}

\nc{\mchar}{\mrm{Char}}
\nc{\Hom}{\mrm{Hom}}
\nc{\id}{\mrm{id}}

\nc{\remark}{\noindent{\bf{Remark:}}}
\nc{\remarks}{\noindent{\bf{Remarks:}}}

 \nc{\delete}[1]{}
 \nc{\grad}[1]{^{({#1})}}
 \nc{\fil}[1]{_{#1}}

\nc{\BA}{{\Bbb A}} \nc{\CC}{{\Bbb C}} \nc{\DD}{{\Bbb D}}
\nc{\EE}{{\Bbb E}} \nc{\FF}{{\Bbb F}} \nc{\GG}{{\Bbb G}}
\nc{\HH}{{\Bbb H}} \nc{\LL}{{\Bbb L}} \nc{\NN}{{\Bbb N}}
\nc{\PP}{{\Bbb P}} \nc{\QQ}{{\Bbb Q}} \nc{\RR}{{\Bbb R}}
\nc{\TT}{{\Bbb T}} \nc{\VV}{{\Bbb V}} \nc{\ZZ}{{\Bbb Z}}
\nc{\Cal}[1]{{\mathcal {#1}}}
\nc{\mop}[1]{\mathop{\hbox {\rm #1} }}
\nc{\smop}[1]{\mathop{\hbox {\eightrm #1} }}
\nc{\mopl}[1]{\mathop{\hbox {\rm #1} }\limits}
\nc{\frakg}{{\frak g}}
\nc{\g}[1]{{\frak {#1}}}
\def \restr#1{\mathstrut_{\textstyle |}\raise-8pt\hbox{$\scriptstyle #1$}}
\def \srestr#1{\mathstrut_{\scriptstyle |}\hbox to
  -1.5pt{}\raise-4pt\hbox{$\scriptscriptstyle #1$}}
\nc{\wt}{\widetilde}
\nc{\wh}{\widehat}
\nc{\un}{\hbox{\bf 1}}
\nc{\redtext}[1]{\textcolor{red}{\tt #1}}
\nc{\bluetext}[1]{\textcolor{blue}{#1}}
\nc{\comment}[1]{[[{\tt {#1}}]] }
\nc{\R}{{\mathbb R}}

\nc\fleche[1]{\mathop{\hbox to #1 mm{\rightarrowfill}}\limits}
\def\semi{\mathrel{\times}\kern -.85pt\joinrel\mathrel{\raise 1.4pt\hbox{${\scriptscriptstyle |}$}}}
\def\astq{{\ast}_q\; \! }

\begin{document}

\title[Twisted dendriform algebras]
      {Twisted dendriform algebras \\ and the pre-Lie Magnus expansion}

\author{Kurusch Ebrahimi-Fard}
\address{Dept. de F\'{i}sica Te\'orica, Facultad de Ciencias, Universidad de Zaragoza, E-50009 Zaragoza, Spain.}
        \email{kef@unizar.es, kurusch.ebrahimi-fard@uha.fr}         
        \urladdr{http://www.th.physik.uni-bonn.de/th/People/fard/}

\author{Dominique Manchon}
\address{Universit\'e Blaise Pascal,
         C.N.R.S.-UMR 6620
         63177 Aubi\`ere, France}
         \email{manchon@math.univ-bpclermont.fr}
         \urladdr{http://math.univ-bpclermont.fr/~manchon/}

\date{March 29, 2010}

\begin{abstract}
In this paper an application of the recently introduced pre-Lie Magnus expansion to Jackson's $q$-integral and $q$-exponentials is presented. Twisted dendriform algebras, which are the natural algebraic framework for Jackson's $q$-analogs, are introduced for that purpose. It is shown how the pre-Lie Magnus expansion is used to solve linear $q$-differential equations. We also briefly outline the theory of linear equations in twisted dendriform algebras.     

\bigskip

\smallskip


\end{abstract}

\maketitle
\tableofcontents

\section{Introduction}
\label{sect:intro}

This article is a continuation of our recent work on the pre-Lie Magnus expansion and linear dendriform equations. In  \cite{EM1} we use Loday's dendriform algebra \cite{Lod1} to describe the natural pre-Lie algebra structure underlying the classical Magnus expansion \cite{Magnus}. This point of view motivated us to explore the solution theory of a particular class of linear dendriform equations  \cite{EM2} appearing in the contexts of different applications, such as for instance perturbative renormalization in quantum field theory. Our results fit into recent developments exploring algebro-combinatorial aspects related to Magnus' work \cite{Chap2,Chap3,Gelfand,Iserles00}. In both references  \cite{EM1,EM2} we mentioned Jackson's $q$-integral and linear $q$-difference equations as a particular setting in which our results apply. 

In the current paper we would like to explore in more detail linear $q$-difference equations and their ($q$-)exponential solutions in terms of a possible $q$-analog of the pre-Lie Magnus expansion. We are led to introduce the notion of twisted dendriform algebra as a natural setting for this work.

\smallskip

The paper is organized as follows. The introduction recalls several
mathematical structures  needed in the sequel such as pre-Lie algebra, bits of
$q$-calculus, and Rota--Baxter algebra. In Section \ref{sect:jackson} linear
$q$-difference equations as well as Jackson's $q$-analog of the Riemann
integral are introduced. The next section contains the definition of unital
twisted dendriform algebra and extends results from earlier work
\cite{EM1,EM2} to this setting. Finally in Section \ref{sect:qmagnus} we
explain how the pre-Lie Magnus expansion in this setting gives rise to a
version of the Magnus expansion involving Jackson integrals. In an appendix we 
briefly remark on a link of our results to finite difference operators.     


\subsection{Preliminaries}
\label{ssect:prelim}

In this paragraph we summarize some well-known facts that will be of use in the rest of the paper.


\subsubsection{Pre-Lie algebras}
\label{ssect:preLie}

Let us start by recalling the notion of pre-Lie algebra \cite{AG, Cartier09, Chap1}. A {\it{left (right) pre-Lie algebra}} is a $k$-vector space $P$ with a bilinear product $\rhd$ ($\lhd$) such that for any $a,b,c \in P$:
 \allowdisplaybreaks{
\begin{eqnarray}
    (a\rhd b)\rhd c - a\rhd (b\rhd c) &=& (b\rhd a)\rhd c - b\rhd(a\rhd c),	\label{prelie1}\\
    (a\lhd b)\lhd c - a\lhd (b\lhd c)  &=& (a\lhd c)\lhd b - a\lhd(c\lhd b).	        \label{prelie2}
\end{eqnarray}}
Observe that for any left pre-Lie product $\rhd$ the product $\lhd$ defined by
$a \rhd b := - b \lhd a$ is right pre-Lie. One shows easily that the left pre-Lie identity rewrites as:
\begin{equation*}
	L_\rhd([a,b])=[L_\rhd (a),L_\rhd (b)],
\end{equation*}
where the left multiplication map $L_\rhd (a) : P \to P$ is defined by $L_\rhd
(a)(b) := a \rhd b$, and where the bracket on the left-hand side is defined by
$[a,b]:=a \rhd b - b \rhd a$. As a consequence this bracket satisfies the
Jacobi identity and hence defines a Lie algebra on $P$, denoted by
$\mathcal{L}_P$: in order to show that (see \cite{AG}) we can add a unit $\un$ to the left
pre-Lie algebra by considering the vector space $\overline P=P\oplus k.\un$
together with the extended product $(a+\alpha\un)\rhd(b+\beta\un)=a\rhd
b+\alpha b+\beta a+\alpha\beta\un)$, which is still left pre-Lie. The map
$a\mapsto L_\rhd(a)$ is then obviously an injective map from $P$ into
$\mop{End}\overline P$, which preserves both brackets.

An easy but important observation is that a commutative left (or right) pre-Lie algebra is necessarily associative, see \cite{AG}.


\subsubsection{Rota--Baxter algebras}
\label{ssect:RBalg}

Recall \cite{Baxter,At,Rota} that a {\it{Rota--Baxter algebra}} (over a field $k$) {\it{of weight}} $\theta \in k$ is an associative $k$-algebra $A$ endowed with a $k$-linear map $R: A \to A$ subject to the following relation:
\begin{equation}
\label{RB}
    R(a)R(b) = R\bigl(R(a)b + aR(b) + \theta ab\bigr).
\end{equation}
The map $R$ is called a {\sl Rota--Baxter operator of weight $\theta$\/}. The map $\widetilde{R}:=-\theta id - R$ also is a weight $\theta$ Rota--Baxter map on $A$. Both, the image of $R$ and $\tilde{R}$ form subalgebras in $A$. Recalling the classical integration by parts rule one realizes that the ordinary Riemann integral, $If(x):=\int_0^xf(y)dy$, is a weight zero Rota--Baxter map. Let $(A,R)$ be a Rota--Baxter algebra of weight $\theta$. Define:
\begin{eqnarray}
\label{RBasso}
             a \ast_R b := R(a)b + aR(b) + \theta ab.
\end{eqnarray}
The vector space underlying $A$, equipped with the product $\ast_R$ is again a Rota--Baxter algebra of weight $\theta$ and $R(a *_R b)=R(a)R(b)$.  Now observe that Rota--Baxter algebras of any other kind than associative still make sense, with the same definition except that the associative product is replaced by a more general
bilinear product. Indeed, let $A$ be an associative Rota--Baxter algebra of weight $\theta$ and define:
\begin{eqnarray}
\label{RBpre-Lie}
            a \rhd_R b := [R(a),b] - \theta ba = R(a)b - b R(a) - \theta ba.
\end{eqnarray}
The vector space underlying $A$, equipped with the product
$\rhd_R$ is a Rota--Baxter left pre-Lie algebra in that sense. Observe that:
$$
	R(a *_R b) + R(b \rhd_R a) = R(R(a)b) + R(R(b)a).
$$


\subsubsection{The pre-Lie Magnus expansion} 
\label{ssect:magnus}

Let $A$ be the algebra of piecewise smooth functions on $\R$ with values in some associative algebra, e.g. square matrices. Recall Magnus' expansion \cite{Magnus}, which allows us to write the formal solution of the initial value problem:
\begin{equation}
\label{ivp}
	\dot{X}(t) = U(t)X(t), 
	\hskip 10mm 
	X(0)=\un,
	\hskip 10mm 
	U \in \lambda A[[\lambda]]
\end{equation}
as $X(t)=\exp(\Omega(U)(t))$. The function $\Omega(U)(t)$ solves the differential equation: 
\begin{equation}
\label{magnus}
	\dot\Omega(U)(t) = \frac{ad_\Omega}{e^{ad_\Omega} -1}(U)(t) 
		                     = U(t) + \sum_{n > 0} \frac{B_n}{n!} ad^n_{\Omega}U(t),
\end{equation}
with initial value $\Omega(U)(0)=0$, where the $B_n$ are the Bernoulli numbers. Here, $e^{ad_\Omega}$ denotes the usual formal exponential operator series. We refer the reader to the recent works \cite{Iserles00,Iserles02,BCOR09} for more details on Magnus' result and its wide spectrum of applications.  In fact, Magnus' expansion allows us to rewrite the {\it{Dyson--Chen series}} of iterated integrals, seen as a formal solution of the initial value problem (\ref{ivp}): 
\allowdisplaybreaks{
\begin{eqnarray*}
	X(t) &=& \un 
			+ IU(t) \
	   		+ I\big(UI(U)\big)(t)
	   		+ I\big(UI(UI(U))\big)(t) + \cdots \\
	       &=&\exp\big(\Omega(U)(t)\big),
\end{eqnarray*}}
where $I$ stands for the Riemann integral $IU(t):=\int_0^tU(s)\,ds$. Let us remark that the sum of iterated commutators on the right-hand side of the second equality in (\ref{magnus}) accounts for the non-commutativity of the underlying function algebra $A$. We now write the Magnus expansion using the natural left pre-Lie product $U \rhd_I V:=[IU,V]$, $U,V \in A$, implied by the weight zero Rota-Baxter relation \eqref{RBpre-Lie}:
\begin{equation}
\label{eq:magnus-bis}
	\dot\Omega(U)(t) = \frac{L_{\rhd_I}(\dot{\Omega})}{e^{L_{\rhd_I}(\dot{\Omega})} -1}(U)(t) 
					    = U(t) +  \sum_{n > 0} \frac{B_n}{n!} L^n_{\rhd_I}({\dot{\Omega}})U(t),
\end{equation} 
where $L_{\rhd_I}(U)V:=U\rhd_I V$. We call this the {\it{pre-Lie Magnus expansion}} \cite{EM1,EM2}, see paragraph \ref{ssect:preLieMag} below. The first few terms are:
\begin{eqnarray*}
\dot{\Omega}(U) &=& U -\frac 12 U\rhd_I U 
					+\frac 14 (U\rhd_I U)\rhd U+\frac 1{12} U\rhd_I(U \rhd_I U) \\
			 & & \qquad\  -\frac 16   ((U\rhd_I U)\rhd U )\rhd_I U - \frac{1}{12} U\rhd_I((U \rhd_I U) \rhd_I U )
					    + \cdots
\end{eqnarray*}

Magnus' expansion is also known as the continuous Baker--Campbell--Hausdorff (\BCH) formula \cite{Wil67, MielPleb,St87,Gelfand} because choosing the function:
$$
	U(t)=\begin{cases}
			x & 0<t<1,\\
			y & 1<t<2,
		\end{cases}		
$$ 
in the above initial value problem gives $\Omega(U)(2)=\BCH (x,y)$, where:
\begin{equation}
\label{clBCH}
	\BCH (x,y) = x + y + \frac 12[x,y] + \frac 1{12}([x,[x,y]]+[y,[y,x]])+\cdots
\end{equation}
stands for the famous {\it{Baker--Campbell--Hausdorff formula}}, appearing in
the product of two exponentials, $\exp(x)\exp(y)=\exp(\BCH (x,y))$. Following
\cite{AG} the pre-Lie Magnus expansion of \cite{EM1} can be approached as
follows: for any pre-Lie algebra $(A,\rhd)$, which is then a Lie algebra, the
\BCH\ formula endows the algebra $\lambda A[[\lambda]]$ with a structure of
pro-unipotent group. This group admits a transparent presentation as follows:
introduce the unit $\un$ again, and define $W: \lambda A[[\lambda]] \to \lambda A[[\lambda]]$ by: 
\begin{equation}
\label{eq:W}
	W(a) := e^{L_\rhd(a)}\un - \un = a + \frac 12 a\rhd a + \frac 16 a\rhd(a\rhd a) + \cdots.
\end{equation}
In fact, one can show that the application $W$ is a bijection, and, in view of \eqref{eq:magnus-bis}, its inverse $W^{-1}$ coincides with $\dot\Omega$ in the case of the pre-Lie product $\rhd_I$ above. We adopt the notation $W^{-1} = \Omega'$ for a general pre-Lie algebra $(A,\rhd)$. Transferring the \BCH\ product by means of the map $W$, namely:
\begin{equation}
\label{diese1}
	a \# b = W \Big(\BCH \big(\Omega'(a),\Omega'(b)\big)\Big),
\end{equation}
we have $W(a) \# W(b) = W\big(\BCH (a,b)\big)=e^{L_\rhd(a)}e^{L_\rhd(b)}\un-\un$ and:
\begin{equation}
\label{BCHpreLie}
	\BCH(\Omega'(a),\Omega'(b))=\Omega'(a \# b).
\end{equation}
Hence $W(a) \# W(b) = W(a) + e^{L_\rhd(a)}W(b)$. The product $\#$ is thus given by the simple formula\footnote{We first learned about this elegant formula and its link to the \BCH\ expansion from private communication with M. Bordemann.}:
\begin{equation}
\label{diese2}
	a\# b=a+e^{L_\rhd(\Omega'(a))}b.
\end{equation}
In particular, when the pre-Lie product $\rhd$ is associative this simplifies to $a\#b = a \rhd b + a + b$. The inverse is given by $a^{\#-1}=W\big(-\Omega'(a)\big)=e^{-L_\rhd(\Omega'(a))}\un - \un$.


\subsubsection{$q$-numbers and $q$-exponentials}
\label{sssect:q-calc}

For details we refer the reader to \cite{KC02}. Recall the definition of $q$-{\it{numbers}}:
$$
	[k]_q:=\frac{1-q^k}{1-q} = \sum_{i=0}^{k-1} q^i,
$$
$0<q<1$ and $k \in \mathbb{Z}_+$, and $q$-{\it{factorials}}:
$$
		[k]_q!:=\begin{cases}
					1 \hbox{ for } k=0,\\
					\prod_{i=1}^k [i]_q \hbox{ for } k\not=0.
			   \end{cases}
$$ 
Jackson's $q$-{\it{exponential functions}} are defined as:
$$
	e_q(x):= \sum_{n = 0}^\infty \frac{x^n}{[n]_q!} 
	\quad {\rm{and}} \quad
	E_q(x):= \sum_{n = 0}^\infty q^{\frac{n(n-1)}{2}} \frac{x^n}{[n]_q!}.  
$$
Observe that $E_q(x)=e_{1/q}(x)$. It has been shown that if the $q$-commutator
$[x,y]_q:=xy-qyx$ vanishes one finds $e_q(x)e_q(y)=e_q(x+y)$. However, if the $q$-commutator of $x$ and $y$ does not vanish, then the product of two $q$-exponentials is encoded by a $q$-analog of the {\BCH}-expansion: 
\begin{equation}
\label{qBCH1}
	e_q(x)e_q(y)=e_q\big(\BCH_q(x,y)\big),
\end{equation}
with:
\begin{equation}
\label{qBCH2}
	\BCH_q(x,y):= x + y - \frac{1}{[2]_q} [y,x]_q - \frac{q}{[2]_q[3]_q!}\big( [[x,y]_q,x]_q +  [y,[x,y]_q]_q \big) + \cdots
\end{equation}
For more details we refer the reader to \cite{KS91,KS94} and \cite{DK95}, where the first four terms of the $q$-\BCH\ expansion have been calculated explicitly. Further below we will show the link to a $q$-analog of the pre-Lie Magnus expansion. For an account of a $q$-analog of the logarithm, see \cite{KvA}.


\section{$q$-differences and Jackson's $q$-integral}
\label{sect:jackson}

Let $A$ be the algebra of continuous functions on ${\mathbb R}$ with values in some not necessarily commutative unital algebra $B$. As a guiding example we have $n \times n$-matrices with complex entries in mind. Let $0<q < 1$, and consider the $q$-difference operator on $A$:
\begin{equation}
\label{qDiff1}
	\partial_qF(t):=\frac{F(qt) - F(t)}{(q-1)t}.
\end{equation}
This operator satisfies the following $q$-twisted Leibniz rule:
\begin{equation}
\label{q-leibniz}
	\partial_q(FG)(t)=\partial_q(F)(t)M_qG(t) + F(t)\partial_q(G)(t),
\end{equation}
where $M_qF(t):=F(qx)$ is the $q$-dilation operator. Observe that the $q$-Leibniz identity is equivalent to: 
\allowdisplaybreaks{
\begin{eqnarray*}
            \partial_{q}(FG)(t)= \partial_{q}(F)(t) G(t) + F(t)\partial_{q}(G)(t) + t(q-1)\partial_{q}F(t) \partial_{q}G(t). 
\end{eqnarray*}} 
Using the $q$-dilation operator $M_q$ we may write:
\begin{equation}
\label{qDiff2}
	\partial_qF(t):=\frac{(id- M_q)F(t) }{(id- M_q)\iota(t)}.
\end{equation}
Here $\iota(t)=t1_B$ denotes the identity map and $1_B$ is the unit in $B$, which we will suppress in the sequel for the sake of notational clarity.

{\exam
\label{ex1d}{
$\partial_qt^k=[k]_qt^{k-1}$, $k \in \mathbb{Z}_+$.
}}

\smallskip

From the foregoing formula (\ref{qDiff2}) we immediately deduce the inverse operation, or anti-derivation, given by Jackson $q$-integral:
\begin{equation}
\label{Jackson}
	I_qf(t) := \int_0^tf(y)\,d_qy 
	             := (1-q) \sum_{k\ge 0}q^ktf(q^kt).
\end{equation}
A direct check shows that $\partial_q \circ I_q(f)=f$, and $I_q \circ \partial_q(f)(t)=f(t)-f(0)$. Setting $F=I_qf$ and $G=I_qg$ in \eqref{q-leibniz} yields the $q$-twisted weight zero Rota--Baxter type relation:
\begin{equation}
\label{qrb}
	I_q\big(f M_qI_q(g)+I_q(f) g\big)=I_q(f) I_q(g).
\end{equation}

Jackson's $q$-integral may be written in a more algebraic way, using again the $q$-dilation operator $M_q$ \cite{At,Rota1}:
\begin{equation}
\label{Poperator}
    P_qf(t) := \sum_{n>0}M_q^{n}f(t) = \sum_{n>0}f(q^{n}t).      
\end{equation}

\begin{prop}
The maps $P_q$ and $id+P_q=:\hat{P}_q$ are Rota--Baxter operators of weight $1$ and $-1$, respectively.
\end{prop}

\begin{proof}
For the map $P_q$ we find
\allowdisplaybreaks{
\begin{eqnarray*}
     P_q(f) P_q(g) 	&=& \sum_{n>0}M_q^{n}(f) \: \sum_{m>0}M_q^{m}(g)  \\
                  	      	&=& \sum_{n,m>0}M_q^{n+m}(f) M_q^{m}(g) 
	      				+ \sum_{m,n>0} M_q^{m}(f) M_q^{n+m}(g) 
						+ \sum_{n>0} M_q^{n}(fg)                              \\
                 	 	&=& P_q(P_q(f)\: g) + P_q(f \:P_qg) + P_q(fg).                                  
\end{eqnarray*}}
The second assertion follows easily by an analogous calculation.  
\end{proof}

The Jackson $q$-integral is given in terms of the above operator as:
\begin{equation}
     I_qf(t) = (1-q) \hat{P}_q (\iota f)(t).
\end{equation}
Recall that $\iota: t \mapsto \iota(t)=t$ is the identity function. One shows that the Jackson $q$-integral satisfies the identity:
\begin{equation*}
    I_q(f)  I_q(g) + (1-q)I_q(f \iota g) = I_q \Big(I_q (f) g + f  I_q(g) \Big).
\end{equation*} 

{\exam
\label{ex1i}{
$ \int_0^t y^k\,d_qy = \frac{t^{k+1}}{[k+1]_q}$, $k\in\mathbb{N}$. 
}}


\subsection{Linear $q$-difference equations}
\label{ssect:qdiffeq}

We briefly recall the well-known results for linear homogenous matrix $q$-difference equations. First, observe that Jackson's $q$-exponentials are eigenfunctions for the $q$-difference operator:
$$
	\partial_q e_q(t)=e_q(t) 	
	\qquad {\rm{and}} \qquad  
	\partial_{\frac 1q} E_{q}(t)=E_{q}(t).
$$
Setting $y(t):=e_q(t)E_q(-t)$ and applying the twisted Leibniz rule we easily compute \cite{B06} $\partial_q y=0$, hence $t\mapsto y(t)$ is a constant function. Setting $t=0$ we thus obtain:
\begin{equation}
\label{expinverse}
	e_q(t)E_q(-t)=e_q(t)e_{q^{-1}}(-t)=1.
\end{equation}
We refer the reader to the standard references mentioned above for more details on the $q$-differences and Jackson's $q$-integral, see also \cite{B06}.

Consider the linear homogenous $q$-difference equations:
\begin{equation}
\label{eq:qdiff1}
	\partial_qX(t) = U(t)M_qX(t) 
	\qquad {\rm{and}} \qquad 	
	\partial_qY(t) = Y(t)V(t), 
\end{equation}
with $X(0)=Y(0)=\un$. It is obvious that the above $q$-difference equations correspond to the $q$-integral equations:    
\begin{equation*}
	X(t) = \un + I_q(UM_qX)(t) 
	\qquad {\rm{and}} \qquad 	
	Y(t) = \un + I_q(YV)(t). 
\end{equation*} 

As a remark we mention the solution in form of an infinite product \cite{B06}:
$$
	X(t)=\prod_{n \ge 0} \big(\un + (1-q)q^ntU(q^nt)\big).
$$
From this one easily verifies that if $U$ and $V$ are constant matrices, the solutions to the above autonomous $q$-difference equations are given by:
\begin{equation}
\label{constsol}
	X(t)=E_q(tU)
		\qquad {\rm{and}} \qquad 
	Y(t)=e_q(tV),
\end{equation}
respectively.

\begin{thm} \label{thm:qDiffFact} 
Let $X,Y$ and $U$ be such that $Y(0)X(0)=\un$ and:
$$
	\partial_q X(t) = U(t)M_qX(t)
		\qquad {\rm{and}} \qquad 
	\partial_q Y(t) = - Y(t)U(t)
$$
then 
$$
	Y(t)X(t)=\un.
$$
\end{thm}
\begin{proof}
A simple direct check shows that the product Z=YX verifies $\partial_qZ=0$, hence $t\mapsto Z(t)$ is a constant. Setting $t=0$ gives the result.
\end{proof}

Applying Theorem \ref{thm:qDiffFact} to a constant matrix $U$ and setting $t=1$ we immediately find the following corollary.

\begin{cor}
For any constant matrix $U$ we have:
\begin{equation*}
	e_q(U)^{-1}=E_q(-U)=e_{\frac 1q}(-U).
\end{equation*}
\end{cor}

At this point a natural question suggests itself: can we express the solutions to the above linear $q$-difference equations as the $q$-exponential of some $q$-analog of the pre-Lie Magnus expansion? In this work we suggest, using the pre-Lie Magnus expansion itself with respect to
Jackson's $q$-integral, to put the solutions in ordinary exponential form.  In fact, our approach in general is based on a formula which can be traced back at least to Hardy and Littlewood (\cite{HL}, see \cite{NG} eq. (53) and \cite{Quesne03} Paragraph 3) which allows us to rewrite $q$-exponentials as ordinary exponentials. However, let us emphasize that it does not suffice to replace the Riemann integral in Magnus' expansion (\ref{magnus}) by Jackson's $q$-integral. Indeed, considering the results presented in this section compared to those in references \cite{EM1,EM2} we are led to the notion of twisted dendriform algebra to be introduced in the next section.


\section{Twisted dendriform algebras}
\label{sect:twistdend}

The notion of dendriform algebras was introduced by Loday~\cite{Lod1}. It developed over the past decade into an active research branch of algebra and combinatorics. The basic idea is to express the multiplication of an associative algebra by a sum of operations together with a set of relations so that the associativity of the multiplication follows from the sum of these relations. As it turns out, a dendriform algebra may at the same time be seen as an associative as well as a pre-Lie algebra. Main examples of dendriform algebras are provided by the shuffle and quasi-shuffle algebra as well as associative Rota--Baxter algebras. However, in the foregoing section we have seen that Jackson's $q$-integral is almost a Rota--Baxter map. This motivates the introduction of a proper extension of Loday's notion, called twisted dendriform algebra.


\subsection{General setting}
\label{ssect:general}

Before giving the general definition we start with an example coming from Jackson's $q$-integral. In view of identity (\ref{qrb}) in the preceding section we are led to introduce the two following bilinear operations on the algebra $A$ of continuous functions:
\begin{equation}
\label{qdend}
	f\succ_q g:=I_q(f) g,\hskip 20mm f\prec_q g:=f I_q(g).
\end{equation}
We also set:
\begin{equation}
\label{dendasso1}
	f\astq g := I_q(f)g + f M_q I_q(g).
\end{equation}
In view of the simple relation $M_q  I_q = q I_q M_q$, this can also be written as:
\begin{equation}
\label{qdendasso2}
	f\astq g:=f\succ_qg + f\prec_q \Phi_q(g),
\end{equation}
with $\Phi_q := qM_q$. From \eqref{qrb} the following three equations can be derived:
\allowdisplaybreaks{
\begin{eqnarray}
	(f\prec_q g)\prec_q h&=&f\prec_q(g\ \astq h),\\
	(f\succ_q g)\prec_q h &=& f\succ_q (g\prec_q h),\\
	f\succ_q(g\succ_q h)&=&(f\ \astq g)\succ_q h.
\end{eqnarray}}
Observe that the product $\astq$ is associative: indeed, thanks to \eqref{qrb} we have:
\begin{eqnarray*}
	I_q\big((f\ \astq g)\ \astq h\big)=I_q(f)I_q(g)I_q(h)=I_q\big(f\ \astq (g\ \astq h)\big).
\end{eqnarray*}
The three relations above are then reminiscent of the axioms of a dendriform algebra \cite{Lod1}, except that the associative product $\astq$ is not the sum of the two operations $\prec_q$ and $\succ_q$. We however remark that $\Phi_q$ is an isomorphism for both of the binary laws $\prec_q$ and $\succ_q$, as well as for the associative product $\astq$. 

Moreover,  one verifies that: 
\begin{equation}
\label{qdendpreLie}
	f \rhd_q g := f\succ_q g - g \prec_q \Phi_q(f)
\end{equation} 
is left pre-Lie. It can be checked directly, but it is a simple consequence of Lemma \ref{lem:preLieDend} below. Observe that:
$$
	I_q(f *_q g + g \rhd_q f) = I_q(f\succ_q g + g \succ_q f) = I_q(I_q(f)g) + I_q(I_q(g)f),
$$
which asks for a generalization of the Bohnenblust-Spitzer identity \cite{EMP09}, \cite{NT09} in the context of Jackson's $q$-integrals. In view of the next section let us do some simple calculations to get familiar with this product.  Hence, let $a,b \in A$ be constant functions. Then:
\allowdisplaybreaks{
\begin{eqnarray}
	a \rhd_q b (t) &=&  a\succ_q b (t)- b \prec_q \Phi_q(a) (t)  \nonumber\\
	                       &=& I_q(a)(t)b - qbI_q(a)(t)  \nonumber\\
	                       &=& abt - qba t=[a,b]_qt, \label{qcalc1}\\
	a \rhd_q (a \rhd_q b)(t) &=& a \succ_q (a \rhd_q b)(t) - (a \rhd_q b) \prec_q \Phi_q(a)(t)  \nonumber\\
						  &=& I_q(a)(t) (a \rhd_q b)(t) - q(a \rhd_q b)(t) I_q (a)(t) \nonumber\\
						  &=& a[a,b]_qt^2 - q[a,b]_qat^2 = [a,[a,b]_q]_qt^2 \label{qcalc1}.
\end{eqnarray}}
As a useful remark, we note that if $a,b$ commute classically, i.e. $ab=ba$, then we find $a \rhd_q b (t) =b \rhd_q a (t)=(1-q)abt$ and:
\allowdisplaybreaks{
\begin{eqnarray*}
	 (a \rhd_q a) \rhd_q b(t) &=& I_q(a \rhd_q a)(t)b - bM_qI_q(a \rhd_q a)(t)\\
							&=& (1-q)a^2b\frac{t^2}{[2]_q} - q^2(1-q)a^2b\frac{t^2}{[2]_q}\\
							&=& (1-q^2)(1-q)a^2b\frac{t^2}{[2]_q}=(1-q)^2a^2bt^2\\
	&=&  [a,[a,b]_q]_qt^2 = a \rhd_q (a \rhd_q b)(t).
\end{eqnarray*}}   

The following two simple lemmas will be of use at the end of the paper. Recall that $A$ is the algebra of continuous functions on ${\mathbb R}$ with values in some not necessarily commutative (in the classical sense) unital algebra $B$. 

\begin{lem}
\label{lem:comprelie}
Let $a_1, \ldots, a_n \in A$, constant functions, all commuting in the classical sense, i.e. $a_na_m=a_ma_n$ for all $n,m$. On the subspace spanned by these elements the pre-Lie product is commutative and hence associative.
\end{lem}

\begin{lem}
\label{lem:jackPoly}
Let $P\in A$ be a polynomial function $P(t):=\sum_{n = 0}^N a_n t^n$ with
coefficients in $B$ which are classically commuting, i.e. $a_la_m=a_ma_l$ for any $l,m$. Then the $n$-fold pre-Lie product is given by:
$$
	 P \rhd_q (P \rhd_q( \cdots  P \rhd_q P)\cdots )(t) =  (1-q)^{n-1} (P(t))^nt^{n-1} = \sum_{m = 0}^{nN} b_{m,n} t^{m+n-1},
$$  
where the coefficients are:
\allowdisplaybreaks{
\begin{eqnarray}
\label{coeff}
	b_{m,n}:= (1-q)^{n-1}\sum_{k_0, k_1,\ldots,k_{N}\ge 0 
					\atop
				 { k_0 + \cdots + k_N=n,\atop k_1+2k_2+\cdots + Nk_N=m}} \left( \begin{matrix}  n \\
				    k_0,\ldots,k_N  \end{matrix}\right) a^{k_0}_{0}\cdots a^{k_m}_{m}
\end{eqnarray}}
\end{lem}

\begin{proof}
The proof follows by induction. First, observe that:
\allowdisplaybreaks{
\begin{eqnarray*}
	  P \rhd_q P(t) &=& I_q(P)(t)P(t)-P(t)M_qI_q(P)(t) \\
	  			    &=& P(t)\big(\sum_{n = 0}^N a_n \frac{t^{n+1}}{[n+1]_q} -  a_n q^{n+1} \frac{t^{n+1}}{[n+1]_q}\big)\\
				    &=& (1-q) (P(t))^2 t
\end{eqnarray*}}   
Hence, in general for a $n+1$-fold product:
\allowdisplaybreaks{
\begin{eqnarray*}
	 P \rhd_q (P \rhd_q( \cdots  P \rhd_q P)\cdots )(t)  &=&  (1-q)^{n-1} (P(t))^nt^{n-1}\big( I_q(P)(t) - M_qI_q(P)(t)\big) \\
				    &=& (1-q)^n (P(t))^{n+1} t^n.
\end{eqnarray*}}   
The coefficients (\ref{coeff}) follow from the multinomial formula. 
\end{proof}

In view of Theorem \ref{thm:qDiffFact} and the subsequent question concerning a $q$-analog of Magnus' expansion, the following remark is in order: indeed, when comparing with (\ref{RBpre-Lie}) we observe that the above pre-Lie product cannot be written using the natural Lie bracket on $A$ and the $q$-integral due to the presence of the isomorphism $\Phi_q$. As a result, even for functions with commuting images, e.g. scalar-valued ones, there is no reason for the pre-Lie product \eqref{qdendpreLie} to vanish unless $q=1$. This leads to the following definition: 

\begin{defn} \label{def:twisted Dend}
A {\rm twisted dendriform algebra\/} on some field $k$ is a $k$-vector space $D$ together with two bilinear operations $\prec$ and $\succ$ and a linear automorphism
$\Phi$ such that for any $x,y,z\in D$:
\allowdisplaybreaks{
\begin{eqnarray}
	\Phi(x\prec y)&=&\Phi(x)\prec\Phi(y)		\label{twisted-denda}\\
	\Phi(x\succ y)&=&\Phi(x)\succ\Phi(y)	\label{twisted-dendb}\\
	(x\prec y)\prec z&=&x\prec(y*z)			\label{twisted-dend1}\\
	(x\succ y)\prec z&=&x\succ(y\prec z)	\label{twisted-dend2}\\
	x\succ(y\succ z)&=&(x*y)\succ z,		\label{twisted-dend3}
\end{eqnarray}}
where:
\begin{equation}
\label{asso}
	x\ast y := x\succ y + x\prec \Phi(y).
\end{equation}
\end{defn}

The usual notion of dendriform algebra is recovered when $\Phi$ is the identity of $D$. The above example $(A,\prec_q,\succ_q)$ is obviously a twisted dendriform algebra, with $\Phi=\Phi_q$.

\begin{lem} \label{lem:assoDend}
The product $*$ in a twisted dendriform algebra $(D,\prec,\succ,\Phi)$ is associative, and $\Phi$ is an automorphism for $*$.
\end{lem}

\begin{proof}
The second assertion is immediate. We then have for any $x,y,z\in D$:
\allowdisplaybreaks{
\begin{eqnarray*}
	(x*y)*z	&=&	(x*y)\succ z+(x*y)\prec \Phi(z)\\
			&=&	(x*y)\succ z +(x\succ y)\prec \Phi(z)+\big(x\prec \Phi(y)\big)\prec \Phi(z)\\
			&=&	(x*y)\succ z+x\succ\big(y\prec\Phi(z)\big)+x\prec\big(\Phi(y)*\Phi(z)\big)\\
			&=&	x\succ(y\succ z)+x\succ\big(y\prec\Phi(z)\big)+x\prec\Phi(y*z)\\
			&=&	x\succ(y*z)+x\prec\Phi(y*z)\\
			&=&	x*(y*z).
\end{eqnarray*}}
\end{proof}

\begin{lem}\label{lem:preLieDend}
The operation $\rhd$ defined by $x\rhd y:=x\succ y - y\prec \Phi(x)$ is {\rm{left pre-Lie}}.
\end{lem}

\begin{proof}
For any $x,y,z\in D$ we have:
\allowdisplaybreaks{
\begin{eqnarray*}
\lefteqn{(x\rhd y)\rhd z-x\rhd(y\rhd z)=}\\
	& (x\succ y)\succ z -\big(y\prec \Phi(x)\big)\succ z - z\prec\Phi(x\succ y)+z\prec\Phi\big(y\prec \Phi(x)\big)\\
	& - x\succ(y\succ z)+(y\succ z)\prec\Phi(x)+x\succ\big(z\prec\Phi(y)\big)-\big(z\prec\Phi(y)\big)\prec\Phi(x)\\
	&=  -\big(x\prec\Phi(y)\big)\succ z-\big(y\prec\Phi(x)\big)\succ z+y\succ
		z\prec\Phi(x)+x\succ z\prec\Phi(y)-z\prec\Phi(x\succ y+y\succ x).
\end{eqnarray*}}
This expression is symmetric in $x$ and $y$, which proves the claim.
\end{proof}

We also consider the right pre-Lie operation $\lhd$ defined by $x\lhd y:=-y\rhd x=x\prec\Phi(y)-y\succ x$. One easily verifies that the associative operation (\ref{asso}) and the pre-Lie operations $\rhd$, $\lhd$ all define the same Lie bracket:
\begin{equation*}
    [\![a,b]\!]:=a*b-b*a=a\rhd b-b\rhd a=a\lhd b-b\lhd a.
\end{equation*}
Let us finally remark that $\Phi$ is an automorphism for both pre-Lie structures $\rhd$ and $\lhd$.


\subsection{Unital twisted dendriform algebras}
\label{ssect:unitalDend}

Let $\overline D:= D \oplus k.\un$ be our twisted dendriform algebra augmented by a unit $\un$:
\begin{equation}
\label{unit-dend}
    a \prec \un := a =: \un \succ a
    \hskip 12mm
    \un \prec a := 0 =: a \succ \un,
\end{equation}
implying that $a*\un=\un*a=a$. We extend the linear isomorphism $\Phi$ to $\overline D$ by setting $\Phi(\un)=\un$. Note that, like in the untwisted case, $\un*\un=\un$, but that $\un \prec \un$ and $\un \succ \un$ are not defined~\cite{R}, \cite{Chap1}.


\subsection{The pre-Lie Magnus expansion}
\label{ssect:preLieMag}

We generalize the results of \cite{EM1} to the twisted case. Let us formally define the exponential and logarithm map in terms of the associative product, $\exp^*(x):=\sum_{n \geq 0} x^{*n}/n!$, $\log^*(\un+x):=-\sum_{n>0}(-1)^nx^{*n}/n$, respectively. This makes sense if the base field is of characteristic zero and if the twisted dendriform algebra $D$ is replaced by $\lambda D[[\lambda]]$, or more generally by a twisted dendriform algebra $D'$ with a complete decreasing filtration compatible with the structure. The exponential (resp. the logarithm) is then a bijection from $D'$ onto $\un + D'\subset \overline {D'}$ (resp. vice-versa). In the following we first give a recursive expression in $\un + \lambda D[[\lambda]]$ for the logarithm of the solutions of the following two equations for a fixed $a \in D$:
\begin{equation}
\label{eq:prelie}
     X = \un + \lambda a \prec \Phi(X),
    \hskip 12mm
     Y = \un - Y \succ \lambda a.
\end{equation}
in $\overline D[[\lambda ]]$, in terms of the left pre-Lie product $\rhd$. Formal solutions to (\ref{eq:prelie}) are given by:
 \allowdisplaybreaks{
\begin{eqnarray*}
    X &=& \un + \lambda a +\lambda^2 a\prec\Phi(a)+\lambda^3 a\prec\Phi(a\prec \Phi(a))+\cdots \\
    Y &=& \un - \lambda a +\lambda^2 a\succ a -\lambda^3 (a\succ a)\succ a +\cdots
\end{eqnarray*}}
Let us introduce the following operators in $(D,\prec,\succ,\Phi)$, where $a$ and $b$ are any elements of $D$:
 \allowdisplaybreaks{
\begin{eqnarray*}
    L_\prec(a)(b)&:=& a \prec b \hskip 8mm L_\succ(a)(b):= a\succ b \hskip 8mm
    R_\prec(a)(b):= b \prec a \hskip 8mm R_\succ(a)(b):= b\succ a\\
    L_\lhd(a)(b)&:=& a \lhd b \hskip 8mm L_\rhd(a)(b) := a \rhd b \hskip 8mm
    R_\lhd(a)(b):=b \lhd a \hskip 8mm R_\rhd(a)(b):= b \rhd a.
\end{eqnarray*}}

\begin{thm} \label{thm:main}
Let $\Omega':=\Omega'(\lambda a)$, $a \in D$, be the element of $\lambda \overline{D}[[\lambda]]$ such that $X=\exp^*(\Omega')$, where $X$ is the solution of the first equation (\ref{eq:prelie}). Then $Y=\exp^*(-\Omega')$ is the solution of the second equation (\ref{eq:prelie}), and $\Omega'$ has the following recursive equation:
 \allowdisplaybreaks{
\begin{eqnarray}
    \Omega'(\lambda a) &=& \frac{R_\lhd(\Omega')}{1-e^{-R_\lhd(\Omega')}}(\lambda a)=
                \sum_{m \ge 0}(-1)^m \frac{B_m}{m!}R^m_\lhd(\Omega')(\lambda a),\label{main1}
\end{eqnarray}}
or alternatively:
 \allowdisplaybreaks{
\begin{eqnarray}\label{eq:plm-ter}
    \Omega'(\lambda a) &=& \frac{L_\rhd(\Omega')}{e^{L_\rhd(\Omega')}-1}(\lambda a)
            =\sum_{m\ge 0}\frac{B_m}{m!}L^m_\rhd(\Omega')(\lambda a),\label{main3}
\end{eqnarray}}
where the $B_l$ are the Bernoulli numbers.
\end{thm}

\begin{proof}
Notice that (\ref{main3}) can be immediately derived from (\ref{main1}) thanks to $L_\rhd(b)=-R_\lhd(b)$ for any $b\in D$. We prove (\ref{main1}), which can be rewritten as:
\begin{equation}
    \lambda a = \frac{1 - e^{-R_\lhd(\Omega')}}{R_\lhd(\Omega')}(\Omega'(\lambda a)).
\label{main5}
\end{equation}
Given such $\Omega':=\Omega'(\lambda a) \in \lambda \overline{D}[[\lambda]]$ we must then prove that $X:=\exp^*(\Omega'(\lambda a))$ is the solution of $X = \un + \lambda a \prec \Phi(X)$, where $a$ is given by (\ref{main5}). Let us first remark that:
\begin{equation}
    R_\lhd(\Omega') = R_\prec(\Phi(\Omega')) - L_\succ(\Omega'),
\end{equation}
and that the two operators $R_\prec(\Phi(\Omega'))$ and $L_\succ(\Omega')$ commute thanks to the twisted dendriform axiom (\ref{twisted-dend2}). We have then, using the five twisted dendriform algebra
axioms:
 \allowdisplaybreaks{
\begin{eqnarray*}
    \lambda a =  \frac{1-e^{-R_\lhd(\Omega')}}{R_\lhd(\Omega')}(\Omega')
      &=& \int_{0}^1 e^{-sR_\lhd(\Omega')}(\Omega')\,ds\\
      &=& \int_{0}^1 e^{sL_\succ(\Omega')}e^{-s(R_\prec(\Phi(\Omega')))}(\Omega')\,ds\\
      &=& \int_{0}^1 \exp^*(s\Omega') \succ \Omega' \prec \exp^*\big(-s\Phi(\Omega')\big)\, ds.
\end{eqnarray*}}
So we get:
 \allowdisplaybreaks{
\begin{eqnarray}
    \lambda a \prec \Phi(X) &=& \int_0^1 \Big(\exp^*(s\Omega') \succ \Omega' \prec \exp^*\big(-s\Phi(\Omega')\big) \Big)\prec \exp^*\big(\Phi(\Omega')\big)\, ds\nonumber\\
               &=& \int_0^1 \exp^*(s\Omega') \succ \Omega' \prec \exp^*\big((1-s)\Phi(\Omega')\big)\, ds\label{eq:important}\\
               &=& \sum_{n\ge 0} \sum_{p+q=n} \Omega'^{*p} \succ \Omega' \prec \Phi(\Omega')^{*q}
                                              \int_0^1 \frac{(1-s)^q s^p}{p!q!}\, ds.\nonumber
\end{eqnarray}}
An iterated integration by parts shows that:
\begin{equation*}
    \int_0^1 (1-s)^q s^p\,ds = \frac{p!q!}{(p+q+1)!},
\end{equation*}
which yields:
\begin{equation*}
    \lambda a \prec \Phi(X) = \sum_{n\ge 0}\frac{1}{(n+1)!}\sum_{p+q=n}\Omega'^{*p}\succ \Omega' \prec \Omega'^{*q}.
\end{equation*}
On the other hand, we have:
\begin{equation}
    X-\un=\exp^*(\Omega')-\un =\sum_{n\ge 0} \frac{1}{(n+1)!} \Omega'^{*n+1}.
\end{equation}
Equality (\ref{main5}) follows then from the identity:
\begin{equation*}
    \sum_{p+q=n}\Omega'^{*p} \succ \Omega' \prec \Phi(\Omega')^{*q} =\Omega'^{*n+1}
\end{equation*}
which is easily shown by induction on $n$. Analogously, one readily verifies that:
 \allowdisplaybreaks{
\begin{eqnarray*}
    Y \succ \lambda a &=& \int_0^1  \exp^*(-\Omega') \succ \Big(\exp^*(s\Omega') \succ \Omega' \prec \exp^*\big(-s\Phi(\Omega')\big) \Big)\, ds\\
               &=& \int_0^1 \exp^*\big((s-1)\Omega'\big) \succ \Omega' \prec \exp^*\big(-s\Phi(\Omega')\big)\, ds\\
               &=& \sum_{n\ge 0} \sum_{p+q=n} (-1)^{(p+q)}\Omega'^{*p} \succ \Omega' \prec \Phi(\Omega')^{*q}
                                              \int_0^1 \frac{(1-s)^q s^p}{p!q!}\, ds\\
                &=& \sum_{n\ge 0}\frac{(-1)^{(p+q)}}{(n+1)!}\sum_{p+q=n}\Omega'^{*p}\succ \Omega' \prec \Phi(\Omega')^{*q}.
\end{eqnarray*}}
\end{proof}

\begin{rmk}
In view of \eqref{eq:plm-ter}, the map $\Omega':\lambda D[[\lambda]]\to\lambda D[[\lambda]]$ is the inverse of the map $W$ introduced in Paragraph \ref{ssect:magnus}, which justifies the notation chosen.
\end{rmk}


\subsection{Perspectives}
First, we remark that the pre-Lie Fer expansion \cite{EM1} also applies in the setting of twisted dendriform algebras giving the logarithm of the solutions of the two equations in (\ref{eq:prelie}) as an infinite product.  Secondly, it would be interesting to adapt the results of \cite{EM2} on linear dendriform equations to the twisted dendriform setting, namely find solutions of second-order equations like:
\begin{equation*}
	Z=\un+\lambda Z\succ a+\lambda^2(Z\succ b)\succ c
\end{equation*}
or higher-order equations:
\begin{equation*}
	Z=\un + \lambda Z\succ a_{11} + \lambda^2(Z\succ a_{12})\succ a_{22}
			+ \cdots + \lambda^N \big(\cdots (Z\succ a_{N1})\succ a_{N2}\cdots \big)\succ a_{NN}
\end{equation*}
in the twisted dendriform algebra $\lambda D[[\lambda]]$.


\section{A Magnus expansion for Jackson integrals}
\label{sect:qmagnus}

Let us return to our function algebra $A$ endowed with the twisted dendriform structure given in terms of Jackson's $q$-integral $I_q$ and the linear isomorphism $\Phi_q=qM_q$. We denote by $\rhd_q$ the corresponding pre-Lie product, and by $W_q,\Omega'_q:A\to A$ the corresponding maps defined in Paragraph \ref{ssect:magnus}. Recall that $\Omega'_q=W_q^{-1}$.

\begin{prop} \label{prop:Wasso}
Let $a \in A$ be the constant function taking the value $t \mapsto a(t):=a$. Then the (non-constant) function $W_q(a)(t)$ is given by:
\begin{equation*}
	W_q(a)(t)=a\left(\frac{\exp\big((1-q)at\big)-1}{(1-q)at}\right).
\end{equation*}
\end{prop}
\begin{proof}
From \eqref{qdendpreLie} we have $a\rhd_qa =(1-q)a^2t$. A straightforward
induction argument (or, alternatively, a direct application of Lemma
\ref{lem:jackPoly} to the constant polynomial $P=a$) shows:
\begin{equation}
	a\rhd_q\big(a\rhd_q\cdots (a\rhd_q a)...)(t)=\big(L^{k+1}_{\rhd_q}(a)\un\big)(t)=a^{k+1}(1-q)^kt^k,
\end{equation}
which, in view of \eqref{eq:W}, proves the proposition.
\end{proof}

We are interested in the (ordinary) logarithm $\Omega_q$ of the solution $\widehat Y$ of the following $q$-integral equation:
\begin{equation}
\label{eq:qmagnus}
	\widehat Y = 1-\lambda I_q(\widehat Y a)
\end{equation}
in $A[[\lambda]]$. We introduce the extra unit $\un$ and consider the unital twisted dendriform algebra $\overline A = A \oplus k\un$ (note that the new dendriform $\un$ is not related to the constant function $1$, which is the unit of the original associative algebra $A$). We extend $I_q$ by setting $I_q(\un):=1$. In view of the pre-Lie Magnus expansion above, it is clear that $\widehat Y = I_q(Y)$ solves \eqref{eq:qmagnus}, where $Y$ is the solution of:
\begin{equation*}
	Y = \un - Y \succ_q \lambda a
\end{equation*}
in $\overline A[[\lambda]]$. From Theorem \ref{thm:main} we deduce that 
\begin{equation}
\label{qexp}
	Y(t)=\exp^{\astq}(\Omega'_q(\lambda a)(t)),
\end{equation}
where the $q$-analog of Magnus' expansion solves:
\begin{equation}
\label{qMagnus}
	\Omega'_q(\lambda a)(t) = \frac{L_{\rhd_q}(\Omega'_q)}{e^{L_{\rhd_q}(\Omega'_q)}-1}(\lambda a)(t)
            		                         =\sum_{m\ge 0}\frac{B_m}{m!}L^m_{\rhd_q}(\Omega'_q)(\lambda a)(t),
\end{equation}
with the pre-Lie product $\rhd_q$ defined in (\ref{qdendpreLie}). Observe that the $q$-analog of Magnus' expansion does not follow from simply replacing the Riemann integral by Jackson's $q$-integral. Instead, at this point it becomes evident that the pre-Lie structure underlying the classical Magnus expansion has to be changed. 

Recall from Theorem \ref{thm:qDiffFact} that $X=Y^{\astq -1}$ solves the linear dendriform equation:
\begin{equation*}
	X = \un + \lambda a\prec \Phi_q(X).
\end{equation*}
In view of $I_q\circ\Phi_q=q^{-1}\Phi_q\circ I_q$, the element $\widehat X=I_q(X)$ solves the $q$-integral equation:
\begin{equation*}
	\widehat X=1+\lambda q^{-1}I_q\big(a\Phi_q(\widehat X)\big).
\end{equation*}

Note that the exponential in (\ref{qexp}) is not the $q$-exponential. Hence the solutions to the $q$-difference respectively $q$-integral equations is given in terms of the classical exponential and all the $q$-deformation structure is encoded in the function:
$$
	\Omega_q(a)(t):=I_q(\Omega'_q(a))(t).
$$
In view of Proposition \ref{prop:Wasso} and equations (\ref{constsol}) as well as of formulas (3.6) and (3.7) in \cite{Quesne03} (cf. Appendix A therein) we have the following results:

\begin{prop}
\label{prop:q-expMagnus1}
Let $a \in A$ be a constant function. Then:
\begin{enumerate}
\item[a)]\label{q-expMagnus2}
\begin{equation*}
	\exp\big(\Omega_q(a)(t)\big) = E_q(at)
		\qquad {\rm{and}} \qquad 
	\exp\big(-\Omega_q(-a)(t)\big) = e_q(at).
\end{equation*}
\item[b)]\label{q-expMagnus1}
\begin{eqnarray*}
\Omega_q(a)(t) &=&\sum_{n>0} \frac{(q-1)^{n-1} }{[n]_q}\frac{t^{n}}{n}a^{n}\\
	\Omega'_q(a)(t) &=& \sum_{n>0} (q-1)^{n-1}\frac{ t^{n-1}}{n}a^{n}=-\frac{1}{(q-1)t}\log\big(1-at(q-1)\big).		
\end{eqnarray*}
\end{enumerate}
\end{prop}

Obviously, it would be desirable to find a link to the $q$-analog of the \BCH\ formula (\ref{qBCH1},\ref{qBCH2}) \cite{KS91}. Recall equation (\ref{BCHpreLie}) which implies:
 \allowdisplaybreaks{
\begin{eqnarray*}
	 e_q\big(\BCH_q(at,bt)\big)&=&e_q(at) e_q(bt) \\ 
				&=&	\exp\big(-\Omega_q(-a)(t)\big) \exp\big(-\Omega_q(-b)(t)\big) \\
		                   &=&	I_q\Big( \exp^{\astq} \big(-\Omega'_q(-a)(t)\big) \astq \exp^{\astq} \big(-\Omega'_q(-b)(t)\big)\Big)\\
		                   &=&  I_q\Big( \exp^{\astq} \big(\BCH_{\astq}(-\Omega'_q(-a)(t),\ -\Omega'_q(-b)(t))\big)\Big)\\
		                   &=&	I_q\Big( \exp^{\astq} \big(-\BCH_{\astq}(\Omega'_q(-b)(t),\ \Omega'_q(-a)(t))\big)\Big)\\
	               		&=& I_q\Big( \exp^{\astq} \big(- \Omega'_q(-b\ \# -a)(t)\big)\Big)\\
	     			&=& \exp \big(- I_q( \Omega'_q(-b\ \# -a))(t)\big)\\
	       			&=& \exp\big(-\Omega_q(-b\ \# -a)(t)\big).	
\end{eqnarray*}}
Here, $\BCH_{\astq}$ denotes the classical \BCH\ expansion (\ref{clBCH}) with respect to the Lie bracket defined in terms of the associative product (\ref{qdendasso2}). We used $\BCH(-b,-a)=-\BCH(a,b)$, as well as equation (\ref{diese2}). However, note that the computation does not lead to $e_q(-b\ \# -a)$. Indeed, $-b\ \# -a$ is not a constant function. Moreover, $-b\ \# -a$ is linear in its second argument. 

The following simple corollary follows from the fact that if $a,b \in A$
commute classically (i.e. $ab = ba$) in $A$ then $a \#b = a \rhd b + a + b$. Recall that $\iota(t)=t$ denotes the identity function.
\begin{cor}
Let $a,b \in A$ be constant functions which commute classically, i.e. $ab =
ba$ in $A$. Then  $a\#b(t) = a \rhd b(t) + a + b = a + b + (1-q)t ab$ and:
$$
	 e_q\big(\BCH_q(at,bt)\big) = e_q(at) e_q(bt) = \exp\Big(-\Omega_q\big(- (a + b - (1-q)ab \iota)\big)(t) \Big).
$$
\end{cor} 

Observe that, using the $q$-Leibniz rule, this follows also naturally from the fact that the function $x(t):=e_q(at) e_q(bt)$ solves the $q$-differential equation:
$$
	\partial_q x(t) = \big(a + b + (q-1)abt\big) x(t).
$$
Finally, note that in the light of Lemma \ref{lem:jackPoly} we find for $P(t):=a + b - (1-q)abt$:
 \allowdisplaybreaks{
\begin{eqnarray*}
	-\Omega_q\big(- (a + b - (1-q)ab \iota)\big)(t) &=& \sum_{n>0} -(q-1)^{n-1} \frac{I_q((P(\iota))^{n}\iota^{n-1})}{n}\\
			&=&  \sum_{n>0} -(q-1)^{n-1} \frac{I_q( \sum_{m = 0}^n b_m \iota^{m+n-1})(t)}{n}\\
	         	&=&  \sum_{n>0} -(q-1)^{n-1} \frac{1}{n}\sum_{m = 0}^n \frac{b_m}{[m+n]_q} t^{m+n}.
\end{eqnarray*}}

We finish with a couple of remarks. First, we have introduced a new functional $q$-analog of the exponential function. Indeed, defining the function: 
$$
	\exp_q := \exp \circ (-\Omega_q) : A \to A
$$
we find for $a,b \in B\subset A$ interpreted as constant functions, $e_q(at) = \exp_q(-a)(t)$, such that:
$$
	\exp_q(-a)\exp_q(-b) =  \exp_q(-b\ \# -a).
$$
Second, as in the classical case, one would expect that the q-\BCH\ formula
should derive from a true $q$-analog of the Magnus expansion, i.e. from
expressing the solutions of $q$-difference equations \eqref{eq:qdiff1} as a
$q$-exponential rather than as an ordinary exponential. A candidate for such a
$q$-analog has recently been investigated by F.~Chapoton \cite{Chap2}, but
this question still remains open.

\section*{Appendix: difference operators}

In this appendix we would to indicate the fine distinction between finite difference calculus and q-difference calculus. 
Indeed, in view of relation (\ref{qrb}) the latter motivated the introduction of the notion of twisted dendriform algebra. Whereas the former, as we will briefly outline here, naturally fits into the Rota--Baxter picture and hence relates to ordinary dendriform algebras \cite{EM1, EM2}..  

Let $A$ be the algebra of piecewise continuous functions on ${\mathbb R}$ with values in some not necessarily commutative unital algebra $B$. Let $h>0$, and consider on $A$ the difference operator:
\begin{equation}
\label{Diff1}
	D_hF(t):=\frac{F(t+h) - F(t)}{h}.
\end{equation}
This operator satisfies the following modified Leibniz rule:
\begin{equation}
\label{h-leibniz}
	D_h(FG)=D_hF.G+F.D_hG+hD_hF.D_hG.
\end{equation}
A right inverse is given by the Riemann summation operator $S_h$ defined by:
\begin{equation}
\label{h-riemann}
S_hf(x)= \begin{cases}
h\sum_{k=1}^{[\frac xh]} f(x-hk) 		& \hbox{if}\ x\ge h, \\ 
0							& \hbox{if}\ 0 \le x<h,\\
-h\sum_{k=0}^{-[\frac xh]-1} f(x+hk) 	& \hbox{if}\ x<0.
\end{cases}
\end{equation}
We indeed have $S_hf(x+h)-S_hf(x)=hf(x)$ for any $x\in\RR$, hence $D_hS_hF=F$. A direct computation shows that for any $F\in A$ and  any $x\in\RR$ we have:
\begin{align}
\label{left-inverse}
	S_hD_hF(x)=F(x)-F\left(x-h[\frac xh]\right).
\end{align}
As a consequence and considering the zero-case in \eqref{h-riemann}, putting $F=S_hf$ and $G=S_hg$ in \eqref{h-leibniz} and applying $S_h$ to both sides we get the weight $h$ Rota--Baxter relation for the operator $S_h$:
\begin{equation}
\label{hRB}
	S_h(f.S_hg+S_hf.g+hfg)=S_hf.S_hg.
\end{equation}
Given some $U\in A$, consider the linear homogeneous equation:
\begin{equation*}
	D_hF=FU
\end{equation*}
with initial condition $F(0)=1$. It is equivalent to the equation:
\begin{equation*}
	F(t)-F\left(t-h[\frac th]\right)=S_h(FU)(t).
\end{equation*}
Restricting ourselves to the functions $F$ which are constant on the interval $[0,h[$ this equation takes the form:
\begin{equation*}
	F=1+S_h(FU)
\end{equation*}
where $1$ denotes the constant function equal to $1$ on $\RR$. Recall (\cite{EM1,EM2}) that $A$ is an ordinary dendriform algebra, with:
\begin{equation*}
	F\succ G=S_h(F)G,\hskip 12mm F\prec G=FS_hG+hFG.
\end{equation*}
Adding formally the unit $\un$ to this dendriform algebra and setting $S_h(\un)=1$, this equation is equivalent to the following equation in $\overline A$ (with $S_hX=F$):
\begin{equation}
\label{dendeq}
	X=\un+X\succ U,
\end{equation}
which have been extensively studied (on a formal level) in \cite{EM1} and \cite{EM2}.


\vspace{1cm}
\subsection*{Acknowledgments}

The first named author is supported by a de la Cierva grant from the Spanish government while being on leave from LMIA, Universit\'e de Haute Alsace, Mulhouse, France. We thank A.~Lundervold for useful comments and proof-reading. 



\begin{thebibliography}{abcdsfgh}
{\small{

\bibitem[AG81]{AG}
	  A.~Agrachev, R.~Gamkrelidze, 
	  {\textsl Chronological algebras and nonstationary vector fields\/}, 
	  J.~Sov.~Math. {\bf{17}} No1, 1650-1675 (1981).

\bibitem[Atk63]{At} 
	F.~V.~Atkinson, 
	{\textsl{Some aspects of Baxter's functional equation}}, 
       	J.~Math.~Anal.~Appl. {\bf 7}, 1-30 (1963).

\bibitem[B06]{B06}
	  G.~Bangerezako, 
	  {\textsl{An introduction to $q$-difference equations\/}}, 
	  preprint, Bujumbura (2006).

\bibitem[Bax60]{Baxter}
    	G.~Baxter,
    	{\textsl{An analytic problem whose solution follows from a simple algebraic identity}},
    	Pacific~J.~Math. {\bf{10}}, 731-742 (1960).

\bibitem[BCOR09]{BCOR09}
	S.~Blanes, F.~Casas, J.A.~Oteo, J.~Ros, 
	{\textsl{Magnus expansion: mathematical study and physical applications}},
	Physics Reports  {\bf{470}}, 151-238 (2009).

\bibitem[Cart72]{Cartier72}
	  P.~Cartier,
	  {\textsl{On the structure of free Baxter algebras}}, 
	  Adv.~Math. {\bf{9}}, 253-265 (1972).
	
\bibitem[Cart09]{Cartier09}
	P.~Cartier,
	{\textsl{Vinberg Algebras and Combinatorics}},
	preprint IHES M/09/34. 

\bibitem[Cha02]{Chap1}
    	F.~Chapoton,
    	{\textsl{Un th\'eor\`eme de Cartier--Milnor--Moore--Quillen pour les alg\`ebres dendriformes et les alg\`ebres braces}},
    	J.~Pure and Appl.~Algebra {\bf{16}}, 1-18 (2002).

\bibitem[Cha08]{Chap2}
    	F.~Chapoton,
    	{\textsl{A rooted-trees $q$-series lifting a one-parameter family of Lie idempotents}},
	Algebra \& Number Theory \textbf{3}, no. 6, 611-636 (2009).

\bibitem[Cha09b]{Chap3}
    	F.~Chapoton,
    	{\textsl{Fractions de Bernoulli-Carlitz et op\'erateurs $q$-Zeta}},
	preprint \texttt{arxiv:0909.1694} (2009).

\bibitem[ChLi01]{ChLi01}
    	F.~Chapoton, M.~Livernet,
    	{\textsl{Pre-Lie algebras and the rooted trees operad}},
    	Int.~Math.~Res.~Not. 2001, 395-408 (2001).
	
\bibitem[ChKa02]{KC02}
	P.~Cheung, V.~Kac,  
	{\textsl{Quantum calculus}}, 
	Springer, New-York (2002).

\bibitem[DK95]{DK95}
	G.~Duchamp, J.~Katriel, 
	\textsl{Ordering relations for $q$-boson operators, continued fraction techniques and the $q$-CBH enigma}, 
	J.~Phys.~A: Math.~Gen.~\textbf{28}, no. 24, 7209-7225, (1995).

\bibitem[EMP09]{EMP09}
	K.~Ebrahimi-Fard, D.~Manchon, F. Patras, 
	\textsl{A Bohnenblust-Spitzer identity for noncommutative Rota--Baxter algebras solves Bogoliubov's counterterm recursion}, 
	J.~Noncommutative~Geometry  \textbf{3}, Issue 2, 181-222 (2009).

\bibitem[EM09a]{EM1}
	K.~Ebrahimi-Fard, D.~Manchon, 
	{\textsl{A Magnus- and Fer-type formula in dendriform algebras\/}}, 
	Found.~Comput.~Math. {\bf{9}},  295-316 (2009).

\bibitem[EM09b]{EM2}
	K.~Ebrahimi-Fard, D.~Manchon, 
	{\textsl{Dendriform equations\/}},
	J.~Algebra \textbf{322},  4053-4079 (2009).

\bibitem[GKLLRT95]{Gelfand}
    	I.~M.~Gelfand, D.~Krob, A.~Lascoux, B.~Leclerc, V.~Retakh, J.-Y.~Thibon,
    	{\textsl{Noncommutative symmetric functions}},
    	Adv.~Math. {\bf{112}}, 218-348 (1995).

\bibitem[HL46]{HL}
	G.~H.~Hardy, J.~E.~Littlewood, 
	{\textsl{Notes on the theory of series (XXIV): a curious power-series}}, 
	Proc. Cambridge Phil. Soc. 42, 85-90 (1946).

\bibitem[Ise00]{Iserles00}
    	A.~Iserles, H.~Z.~Munthe-Kaas, S.~P.~N{\o}rsett, A.~Zanna,
    	{\textsl{Lie-group methods}},
    	Acta Numerica {\textbf{9}}, 215-365 (2000).

\bibitem[Ise02]{Iserles02}
    	A.~Iserles,
    	{\textsl{Expansions that grow on trees}},
    	Notices of the AMS {\bf{49}}, 430-440 (2002).

\bibitem[KaSo91]{KS91}
	J.~Katriel, A.~I.~Solomon, 
	{\textsl{A $q$-analogue of the Baker--Campbell--Hausdorff expansion}}, 
	J.~Phys.~A: Math.~Gen. {\bf{24}}, L1139-L1142  (1991).

\bibitem[KaSo94]{KS94}
	J.~Katriel, A.~I.~Solomon, 
	{\textsl{A no-go theorem for a Lie-consistent $q$-Campbell--Baker--Hausdorff expansion}}, 
	J.~Math.~Phys. {\bf{35}} (11), 6172-6178 (1994).

\bibitem[KvA09]{KvA}
	E.~Koelink, W.~Van Assche, 
	\textsl{Leonhard Euler and a $q$-analogue of the logarithm},  
	Proc.~Amer.~Math.~Soc. \textbf{137}, no. 5, 1663-1676 (2009).

\bibitem[Lod01]{Lod1}
    	J.-L.~Loday,
    	{\textsl{Dialgebras}},
    	Lect.~Notes~Math. 1763, Springer, Berlin pp. 7-66 (2001).

\bibitem[Mag54]{Magnus}
    	W.~Magnus,
    	{\textsl{On the exponential solution of differential equations for a linear operator}},
    	Commun.~Pure Appl.~Math. {\bf{7}}, 649-673 (1954).

\bibitem[MP70]{MielPleb}
    	B.~Mielnik, J.~Pleba\'nski,
    	{\textsl{Combinatorial approach to Baker--Campbell--Hausdorff exponents}},
    	Ann.~Inst.~Henri Poincar\'e {\bf{A XII}}, 215-254 (1970).

\bibitem[NG96]{NG}
	C.~A.~Nelson, M.~G.~Gartley, 
	{\textsl{On the two q-analogue logarithmic functions: $\mop{ln}_q(w)$, $\mop{ln}\{e_q(z)\}$}}, 
	J.~Phys.~A: Math.~Gen. \textbf{29}, no. 24, 8099-8115 (1996).

\bibitem[NT09]{NT09}
	J-C.~Novelli, J-Y.~Thibon, 
	\textsl{A one-parameter family of dendriform identities}, 
	J.~Combin.~Theory Ser.~A \textbf{116}, 864-874 (2009).

\bibitem[Qu03]{Quesne03}
	C.~Quesne,
	{\textsl{Disentangling $q$-Exponentials: A General Approach}},
	Int.~J.~Theoretical Phy. {\bf{43}}, 545-559 (2004).

\bibitem[Ron00]{R}
    	M.~Ronco,
    	{\textsl{Primitive elements in a free dendriform algebra}},
    	Contemp.~Math.~{\bf{207}}, 245-263 (2000).
		
\bibitem[Rot69]{Rota}
    	G.-C.~Rota,
    	{\textsl{Baxter algebras and combinatorial identities. I, II.}},
    	Bull.~Amer.~Math.~Soc. {\bf{75}}, 325-329 (1969); \textsl{ibidem}~330-334.

 \bibitem[Rot95]{Rota1} 
 	  G.-C.~Rota,
  	{\textsl{Baxter operators, an introduction}}, 
	  in ''Gian-Carlo Rota on combinatorics'', Contemp. Mathematicians, Birkh\"auser Boston, Boston, MA,  504-512 (1995).
		
\bibitem[Str87]{St87}
	R.~Strichartz,
	{\textsl{The Campbell--Baker--Hausdorff Dynkin formula and solutions of differential equations}}, 
	J.~Funct.~Analysis  {\bf{72}}, 320-345 (1987). 
	
\bibitem[Wil67]{Wil67}
	R.~M.~Wilcox, 
	{\textsl{Exponential operators and parameter differential in quantum physics}}, 
	J.~Math.~Phys. {\bf{8}}, 962-982 (1967).

}}
\end{thebibliography}
\end{document}